\documentclass[11pt]{amsart}
\usepackage{amssymb}
\usepackage{amsmath}
\usepackage{amsfonts}
\theoremstyle{plain}
\newtheorem{thm}{Theorem}[section]

\newtheorem{lmm}[thm]{Lemma}
\newtheorem{cor}[thm]{Corollary}
\newtheorem{dfn}[thm]{Definition}
\newtheorem{rmk}[thm]{Remark}

\font\fiverm=cmr5

\input prepictex 
\input pictex
\input postpictex
 
\def\al{\alpha}
\def\be{\beta}
\def\ga{\gamma}

\def\vep{\varepsilon}

\def\om{\omega}

\def\ov{\overline}

\def\sm{\setminus}

\begin{document}
                  
\bigskip

\title
{A construction of noncontractible simply connected cell-like
two dimensional Peano continua}   

\author{Katsuya Eda}
\address{School of Science and Engineering,
Waseda University, Tokyo 169-8555, Japan}
\email{eda@logic.info.waseda.ac.jp}

\author{Umed H. Karimov}
\address{Institute of Mathematics,
Academy of Sciences of Tajikistan,
Ul. Ainy $299^A$, Dushanbe 734063, Tajikistan}
\email{umed-karimov@mail.ru}

\author{Du\v san Repov\v s }
\address{Institute of Mathematics, Physics and Mechanics,
and Faculty of Education, 
University of Ljubljana, P.O.Box 2964,
Ljubljana 1001, Slovenia}
\email{dusan.repovs@guest.arnes.si}




\begin{abstract}
Using the topologist sine curve we present a new functorial construction
 of cone-like spaces, starting in the category of all path-connected
 topological spaces with a base point and continuous maps, and ending in
 the subcategory of all simply connected spaces.
If one starts by a noncontractible $n$-dimensional Peano continuum for
any $n>0$, then our construction yields a simply connected noncontractible 
$(n\! +\! 1)$-dimensional cell-like Peano continuum. In particular,
 starting with the circle $\mathbb{S}^1$, one gets a noncontractible simply
 connected cell-like 2-dimensional Peano continuum. 
\end{abstract}

\maketitle
\pagestyle{myheadings} \markright{NONCONTRACTIBLE SIMPLY CONNECTED  CELL-LIKE
CONTINUA}

\section{Introduction}

It is well known that all cell-like
polyhedra are contractible. Griffiths \cite{Gr}
constructed a 2-dimensional nonsimply connected
cell-like Peano continuum: Let $\mathbb{H}_1$ 
be the 1-dimensional {\it Hawaiian earrings} with the base point
$\theta$ at which $\mathbb{H}_1$ is not locally simply connected. 
Let $Y = C(\mathbb{H}_1)$ be the cone over $\mathbb{H}_1$. Then
$\mathbb{H}_1$ can be considered as the base of the cone 
$C(\mathbb{H}_1)$ and 
$\theta$ as its base point. The Griffiths space
is then defined as
the bouquet of two copies of $Y$ with respect 
to the point $\theta$.

A generalization of the Griffiths example is analogous --
instead of the 1-dimensional one considers the 2-dimensional 
Hawaiian earrings \cite{E},
i.e. the subspace $\mathbb{H}_2$ of the 
3-dimensional Euclidean space,
$\mathbb{H}_2 = \{ (x_0, x_1, x_2) \in 
\mathbb{R}^3 |\ (x_0 - 1/k)^2 + x_1^2 + x_2^2 = 
(1/k)^2,\ k\in \mathbb{N}\}$ . 
It is easy to see that this generalization of the Griffiths example is a
3-dimensional noncontractible simply connected cell-like Peano continuum. 

The purpose of the present paper is to construct a 
functor $SC(-, -)$ 
from the category of all path connected spaces 
with a base point and continuous
mappings, to the subcategory of all simply 
connected spaces with a base point. The following are our main results:

\begin{thm}\label{thm:main1}
For every path connected space $Z$ with $z_0\in Z$, the space $SC(Z,z_0)$ is
simply-connected.
\end{thm}

\begin{thm}\label{thm:main2}
For every noncontractible space $Z$ with $z_0\in Z$, the space
 $SC(Z,z_0)$ is noncontractible.
\end{thm}

If $Z$ is a Peano continuum, then $SC(Z, z_0)$
is also a Peano continuum.
If $Z$ is an $n$-dimensional metrizable space for $n>0$, then the space
$SC(Z,z_0)$ is $(n+1)$-dimensional. 
If $Z$ is a compact, then $SC(Z, z_0)$ is a compact
space with trivial shape.  

In particular, when $Z$ is the circle $\mathbb{S}^1$, we get the
following: 
\begin{cor}\label{cor:main1}
For any point $z_0$ of the circle $\mathbb{S}^1$, the space
 $SC(\mathbb{S}^1,z_0)$ is a noncontractible simply
 connected cell-like 2-dimensional Peano continuum.
\end{cor}
As a general reference for algebraic topology we refer the reader to
\cite{Spanier:algtop}. 

\section{Preliminaries}
For any two points $A$ and $B$ in the plane $\mathbb{R}^2$, 
$[A,B]$ denotes the linear segment connecting these points. 
For $a,b \in \mathbb{R}$ with $a < b$, $[a,b]$ denotes the closed
interval, $(a, b)$ denotes the open interval and 
$[a, b)$ and $(a, b]$ denote the half-open intervals, as usual. The unit
interval $[0,1]$ will be denoted by $\mathbb{I}$. To avoid 
confusion between an open
interval and an element of the square $\mathbb{I}\times \mathbb{I}$, 
we shall
write $(a; b)$ for the latter, where $a,b\in \mathbb{I}$. 
Our construction is based on the
piecewise linear topologist sine curve $T$ in the plane. Let
$A_n = (1/n;0),$ $B_n = (1/n; 1)$, for $n\in \mathbb{N} = \{1, 2, 3,
\dots \}$, $A = (0;0),$ $B = (0;1)$ be the 
points of the plane $\mathbb{R}^2$. Let $L_{2n-1} = [A_n, B_n]$
and $L_{2n} = [B_n, A_{n+1}]$. The space $T$ is the subspace of
$\mathbb{I}^2$ defined as the union of all segments $L_n$ and $L =
[A,B]$.

Let $Z$ be any space with a base point $z_0$. 
Then the base set of $SC(Z, z_0)$ is the quotient set of $ T\times Z\ \cup\
\mathbb{I}^2 $ obtained by the identification of the points
$(s,z_0)\in T\times Z$ with $s\in T\subset \mathbb{I}^2$ and by
the identification of each set $\{s\}\times Z$ with the one-point set
$\{s\}$ if $s\in L$. 
There is a natural projection $p: SC(Z, z_0)\to \mathbb{I}^2$. To
$p$ there corresponds a pair of functions $p_1$ and $p_2$ such
that $p(z)= (p_1(z); p_2(z))$. 
For $a = (x; y)\in T$ with $x>0$, the set $p^{-1}(a)$ is denoted by
$Z_a$, which is homeomorphic to $Z$, and for $y\in \mathbb{I}$ the set
$p^{-1}_2(\{ y\})$ is denoted by $M_y$. 
Let $O_\vep (a) = p^{-1}(U_\vep (a))$, where  $U_\vep (a)$ is
the open $\vep$-ball with the center at $a \in \mathbb{I}\times
\mathbb{I}$ with respect to the standard metric. 

The topology of $SC(Z, z_0)$ coincides with the quotient topology 
at each point outside
$L$. A basic neighborhood of a point $a=(0; y)\in L$
is of the form $O_\vep (a)$. Therefore, the topology of $SC(Z, z_0)$ is
the quotient topology when $Z$ is compact.  

Obviously, $SC(-,\ -)$ is a functor from the
category of topological spaces with a base point to itself. The
space $SC(Z, z_0)$ is path-connected, path-connected and locally connected,
finite-dimensional, metrizable or compact if $Z$ is path-connected,
path-connected and locally
connected, metrizable or compact, respectively. In
particular, $SC(Z, z_0)$ is a Peano continuum if $Z$ is a Peano
continuum. 

If $Z$ is compact, $SC(Z)$ is a quotient space of $T\times Z\cup
\mathbb{I}^2$ and hence $SC(Z)$ is also compact. 
Next we show that the shape type of $SC(Z)$ is that of the one-point
space, when $Z$ is compact. To see this let $\mathcal{U}$ be
an open cover of $SC(Z)$. By the compactness of $\mathbb{I}$ we have 
$\vep >0$ such that $p^{-1}([0,\vep )\times [a,b])$ is contained in an
element of $\mathcal{U}$ for every $0\le a < b\le 1$ with $b-a < \vep$. 
By the compactness of $Z$ we also have a cover $O_1, \cdots ,O_m$ of $X$
and points $P_1 = A_1, \cdots, P_n$ on $T$ such that
$[P_i,P_{i+1}]\subset T\cap [\vep ,1]\times \mathbb{I}$ and each 
$[P_i,P_{i+1}]\times O_j$ is contained in an element of $\mathcal{U}$. 
Hence we have a refinement of $\mathcal{U}$ whose nerve is
contractible. This yields the conclusion. 

Let $Z$ be an $n$-dimensional metrizable space for $n>0$. 
Then, since $SC(Z) = \mathbb{I}^2 \cup \bigcup _{n=1}^\infty
p^{-1}([1/n,1]\times \mathbb{I})$ and $p^{-1}([1/n,1]\times
\mathbb{I}\cap T)$ is homeomorphic to $\mathbb{I}\times Z$, 
the dimension of $SC(Z)$ is $n+1$ by 
\cite[Theorems 4.1.3 and 4.1.9]{Engelking:dimension} and 
\cite[p.221]{Morita:dimension}. 
Hence $SC(Z)$ is a cell-like, $(n\! +\! 1)$-dimensional
compact metrizable space, if $Z$ is a $n$-dimensional
compact metrizable space \cite{Lacher, MitchellRepovs}.   

A {\it path} in
$X$ is a continuous mapping of the closed interval $[a, b]$
to $X$. We say that two paths are {\it homotopic} if they are
defined on the same domain and are homotopic relative to their
ends.
The {\it composition} of two paths
$f:[a, b]\to X$ and $g: [b,c]\to X$ such that $f(b) = g(b)$ is a path
$h:[a,c]\to X$ which is defined as follows: 
\[
h(t) = \left\{
\begin{array}{ll}
f(t) \quad &\mbox{if } a\le t\le b \\
g(t) \quad &\mbox{if } b\le t\le c.
\end{array}\right.
\]

Let $f:[a, b]\to X$ and $g: [c,d]\to X$ be paths. We say 
that $f$ is equivalent to $g$ and write $f\cong g$
when $f(a+(b-a)t) = g(c+(d-c)t)$ for each $t\in \mathbb{I}$ and define
$\ov{f}$ as $\ov{f}(t) = f(a+b-t)$ for each $a\le t\le b $. 

A {\it loop} with the {\it base point} $x_0$ in a space $X$ is a path 
$f:[a,b]\to X$ for which $f(a)= f(b)= x_0$.

The product of two loops is
defined in the standard way. The constant mapping to
$\{ x_0\}$ is denoted by $c_{x_0}$. 
Let $f:[a, b]\to X$ be a path, $c$ any point in $[a, b]$ and $\al$ any
loop with the base point at $f(c)$. The {\it modification} of the path
along a loop $\al :\mathbb{I}\to X$ is the path $g:[a, b]\to X$ which is
defined for an interval $[t_1, t_2]\subset (a, b)$ such that $c\in [t_1,
t_2]$ as follows: 
\[
g(s) = \left\{
\begin{array}{ll}
f\left((s-a)(c-a)/(t_1-a) + a\right)\quad   &\mbox{if } a\le s\le t_1 \\
\\
\al((s-t_1)/(t_2-t_1)) \quad &\mbox{if } t_1\le s\le t_2\\
\\
f\left((b-s)(b-c)/(b-t_2) + b\right) \quad  &\mbox{if } t_2\le
s\le b.
\end{array}\right.
\]
The definition of the modification of paths depends on the interval 
$[t_1, t_2]$, however all such paths are homotopic.
(For simplicity of the definition we suppose that the domain of a loop
$\al$ is $\mathbb{I}$, but we shall use a variant of the modification
for loops with arbitrary domains in the sequel.)

A homotopy connecting an injective mapping with the constant one is called
a ${\it contraction}$. 
Whenever possible we shall use the symbol $SC(Z)$
instead of $SC(Z, z_0)$.

\section{Proof of Theorem 1.1}

\begin{lmm}\label{lmm:strngretract}
Let $A$ be a strong deformation retract of $X$ and let $\al :
[0, 1] \to X$ be a path with the end points $\al (0)$ and $\al
(1)$ in $A$. Then there exists a path $\al' : [0, 1] \to A\subset
X$ which is homotopic to $\al$.
\end{lmm}

\begin{proof}
The assertion of the lemma follows directly from the definition of
the strong deformation 
retraction.
\end{proof}

\begin{lmm}\label{lmm:smallmove1}
Let $X$ be any space and $\al = (\al_1, \al_2)$ any path in $X\times \mathbb{I}$ with
the end points $\al (0) \in X\times
\{0\}$ and $\al (1)\in X\times \{1\}$. Then there exists a 
path $\al'$ in $X\times  \mathbb{I}$ homotopic to $\al$ and such
that ${\rm Im}(\al')\subset \{\al_1(0)\}\times \mathbb{I} \cup
X\times \{1 \}$.
\end{lmm}

\begin{proof}
Let $H: \mathbb{I}\times \mathbb{I} \to X\times \mathbb{I}$ be the
homotopy which is defined by the following formulas:
\[
H(s, t) = \left\{
\begin{array}{ll}
(\al_1(0), 2s)  &\mbox{if } 0\le s\le t/2 \\
\displaystyle{(\al_1 (\frac{2s-t}{2-t}) ,\
(1-t)\al_2( \frac{2s-t}{2-t}) + t )}\quad &\mbox{if } t/2 
\le s\le 1.
\end{array}\right.
\]

Obviously, $H(0,t) = \al(0)$, $H(1,t) = \al (1)$, $H(s, 0) = \al (s)$
 and 
${\rm Im} (H(- ,
1))\subset \{\al_1(0)\}\times \mathbb{I}\cup X\times \{1\}$  so
$H$ is the desired homotopy connecting $\al$ and $\al' = H(-
, 1).$
\end{proof}

The path $\al'|_{ [0,\ 1/2]}$ is called the
{\it linear\ part} and $\al'|_{ [1/2,\ 1]}$ is called the {\it
residual part} of the path $\al'$.

In the following lemmata we use the symbols of $B_n$, $Z_{B_n}$,
$U_{\vep}(B_n)$, $O_{\delta}(B_n)$ and $M_0$ which were defined in
Section 2.
\begin{lmm}\label{lmm:smallmove2}
Let $f:\mathbb{I} \to SC(Z)$ be any path. Then for every $n\in
\mathbb{N}$ and every $\vep > 0$ there exist a path
$f_{n,\vep}:\mathbb{I}\to SC(Z)$ and a homotopy $H_{n,\vep}:
\mathbb{I}^2\to SC(Z)$ such that:
\begin{itemize}
\item[(1)] $H_{n,\vep}(s, 0) = f(s), \ \  H_{n,\vep}(s, 1) =
f_{n,\vep}(s);$
\item[(2)] ${\rm Im}(f_{n, \vep})\cap Z_{B_n} = \emptyset;$ and
\item[(3)] $H_{n,\vep}(s, t) = f(s)$ if $f(s)\notin O_\vep (B_n).$
\end{itemize}
\end{lmm}

\begin{proof}
Let $\delta$ be a number such that $0<
\delta < \vep$ and $U_{\delta}(B_n)\cap T =
U_{\delta}(B_n)\cap (L_{2n-1}\cup L_{2n}).$
Since $f^{-1}(Z_{B_n})$ is a compact subset of
$f^{-1}(O_{\delta}(B_n))$, there exists a finite 
set of pairwise disjoint intervals 
$\{[a_k, b_k] :  k\in K_n\}$ which cover
$f^{-1}(Z_{B_n})$ in $f^{-1}(O_{\delta}(B_n))$ and whose end
points lie outside $f^{-1}(Z_{B_n})$. Using the modifications of
paths along loops we can assume without loss of generality that
end points of all paths lie on $T$.
For a given $k\in K_n$ consider the path $f|_{[a_k, b_k]}:[a_k, b_k]\to
O_{\delta}(B_n)$. Since $p^{-1}(U_{\delta}(B_n)\cap T)$ is a
strong deformation retract of $O_{\delta}(B_n)$, the path $f|_{[a_k, b_k]}$
is homotopic to a path $f_{n,k}:[a_k, b_k]\to p^{-1}(U_{\delta}(B_n)
\cap T )\subset O_{\delta}(B_n)$ due to Lemma
\ref{lmm:strngretract}.

The space $p^{-1}(U_{\delta}(B_n)\cap T)$ is naturally homeomorphic
to the product of the interval and the space $Z$. 
The product $([f(a_k), B_n]\cup [B_n,f(b_k)])\times Z$ is a
strong deformation retract of $p^{-1}(U_{\delta}(B_n)\cap T)$.
Therefore the path $f:[a_k, b_k]\to O_{\delta}(B_n)$ is
homotopic to a path in $([f(a_k), B_n]\cup [B_n,f(b_k)])\times Z$, again
 by Lemma \ref{lmm:strngretract}.
By Lemma \ref{lmm:smallmove1}, the path $f:[a_k, b_k]\to
O_{\delta}(B_n)$ is homotopic to a path the linear part of
which lies in $L_{2n-1}\cup L_{2n}$ and the residual part of
which does not intersect $Z_{B_n}$. The linear part can be
slightly deformed in $\mathbb{I}\times \mathbb{I}$ to a segment
 $[f(a_k),f(b_k)]$ with fixed ends $f(a_k)$ and $f(b_k)$,
 which does not contain the point $B_n$.
Since the index $k$ is arbitrary and the number of the intervals 
$\{[a_k, b_k]: k\in K_n\}$ is finite we get the desired mapping
$f_{n,\vep}$.
\end{proof}

Next lemma is a direct consequence of Lemma
\ref{lmm:smallmove2}:

\begin{lmm}\label{lmm:smallmove3}
Any loop in $SC(Z)$ with the base point in $M_0$ is homotopic to a loop
in $SC(Z)\sm \bigcup _{n\in \mathbb{N}}Z_{B_n}$.
\end{lmm}

\begin{center}
\beginpicture 
\setcoordinatesystem units <1.000mm,1.000mm>
\setplotarea x from -10 to 110, y from -10 to 110
\setsolid
\setplotsymbol ({\fiverm.})
\plot 0.000 60.000 100.000 60.000 /
\plot 0.000 90.000 100.000 90.000 /
\plot 0.000 20.000 100.000 20.000 /
\plot 25.000 0.000 25.000 100.000 /
\plot 50.000 100.000 25.000 0.000 /
\plot 50.000 0.000 50.000 100.000 /
\plot 100.000 100.000 50.000 0.000 /
\plot 0.000 100.000 0.000 0.000 /
\plot 100.000 100.000 0.000 100.000 /
\plot 100.000 0.000 100.000 100.000 /
\plot 0.000 0.000 100.000 0.000 /
\put {$A_{3,1}$} at 20.000 17.000
\put {$B_{3,0}$} at 20.000 93.000
\put {$A_{3,0}$} at 33.000 17.000
\put {$C_{4}$} at 35.000 57.000
\put {$C_{5}$} at 20.000 57.000
\put {$A_{2,1}$} at 45.000 17.000
\put {$A_{2,0}$} at 64.000 17.000
\put {$C_{3}$} at 55.000 57.000
\put {$C_{2}$} at 83.000 57.000
\put {$B_{2,1}$} at 43.000 93.000
\put {$B_{2,0}$} at 55.000 93.000
\put {$B_{1,1}$} at 92.000 93.000
\put {$A_{1,1}$} at 105.000 17.000
\put {$C_{1}$} at 105.000 57.000
\put {$B_{1,0}$} at 105.000 93.000
\put {$A$} at -5.000 -3.000
\put {$C$} at -5.000 57.000
\put {$B$} at -5.000 103.000
\put {$B_1$} at 105.000 103.000
\put {$B_2$} at 50.000 103.000
\put {$B_3$} at 25.000 103.000
\put {$A_1$} at 105.000 -3.000
\put {$A_2$} at 50.000 -3.000
\put {$A_3$} at 25.000 -3.000
\put {$\bullet$} at 100.000 0.000
\put {$\bullet$} at 50.000 0.000
\put {$\bullet$} at 25.000 0.000
\put {$\bullet$} at 25.000 100.000
\put {$\bullet$} at 50.000 100.000
\put {$\bullet$} at 100.000 100.000
\put {$\bullet$} at 100.000 90.000
\put {$\bullet$} at 100.000 60.000
\put {$\bullet$} at 100.000 20.000
\put {$\bullet$} at 60.000 20.000
\put {$\bullet$} at 80.000 60.000
\put {$\bullet$} at 95.000 90.000
\put {$\bullet$} at 50.000 90.000
\put {$\bullet$} at 50.000 60.000
\put {$\bullet$} at 50.000 20.000
\put {$\bullet$} at 30.000 20.000
\put {$\bullet$} at 40.000 60.000
\put {$\bullet$} at 47.500 90.000
\put {$\bullet$} at 25.000 90.000
\put {$\bullet$} at 25.000 60.000
\put {$\bullet$} at 25.000 20.000
\put {$\bullet$} at 0.000 20.000
\put {$\bullet$} at 0.000 60.000
\put {$\bullet$} at 0.000 100.000
\put {$\bullet$} at 0.000 90.000
\put {$\bullet$} at 0.000 0.000
\endpicture

Figure 1
\end{center}

\begin{lmm}\label{lmm:retract}
$M_0$ is a strong deformation retract of $SC(Z)\sm \bigcup _{n\in
\mathbb{N}} Z_{B_n}$.
\end{lmm}

\begin{proof}
The deformation $D:(SC(Z)\sm \bigcup _{n\in
\mathbb{N}}Z_{B_n})\times \mathbb{I}\to SC(Z)\sm \bigcup _{n\in
\mathbb{N}}Z_{B_n}$ is given by the piecewise linear mapping (linear
 over every triangle $A_nB_nA_{n+1}$ and $A_{n+1}B_{n+1}B_{n})$ which
maps $[A_n,\ B_n]\setminus \{ B_n\}$ and $[A_{n+1},\ B_n]\setminus 
\{ B_n\}$ to the points $A_n$ and $A_{n+1}$, respectively (see Figure 1).
Since the spaces $Z_{B_n}$ have been deleted, $D$ is well-defined
and continuous.
\end{proof}

The following follows from Lemmata~\ref{lmm:smallmove3} and
\ref{lmm:retract} :

\begin{lmm}\label{lmm:canonical}
Let $f$ be a loop in $SC(Z)$ whose base point is in $M_{0}$. Then
$f$ is homotopic to a loop in $M_0$.
\end{lmm}

Before we show the simple connectivity of $SC(Z)$, we 
exhibit a homotopy from the canonical winding to the constant, in case
when $Z$ is the circle in Figure 2. In the remaining part of this
section we shall use the word ``homotopic'' for loops in a weaker sense, that
is, two loops $f$ and $g$ will be considered to be homotopic if there
exists a homotopy $H(-,t)$ such that $H(-,0)=f$, $H(-,1)=g$ and $H(-,t)$
is a loop for each $t$. 
>From a loop as I in Figure 2 we pull the bottom of the loop to the
left. This is the procedure I. Then we pull up to the loop as III
through the one as II, this is the procedure II. Now the loop is in the upper edge without tangles as III. We contract the
loop to the point $B$, this is the procedure III. 
To generalize these simple procedures
I, II, III, we need to describe them more precisely. 

For a loop $\al$ in $Z$ with the base point $z_0$ and a point $u\in T$, let
$\al _u$ be a loop in $Z_u$ induced naturally by the homeomorphism
between $Z$ and $Z_u$, i.e. $\al _u (t) = (u,\al (t))$, and in
particular the base point of $\al _u$ being $u$. 

We call $\be : [a,b]\to SC(Z)$ a {\it basic loop} at $A_n$, if 
there exists a loop $\al$ in $Z$ with the base point $z_0$ such that 
\begin{itemize}
\item[(a)] $\be (a) = \be (b) = A$, $\be ((2a+b)/3) = \be ((a+b)/2) = A_n$; 
\item[(b)] $\be |_{[a,(2a+b)/3]}$ and $\be |_{[(a+b)/2,b]}$ are linear
           mappings; and 
\item[(c)] $\be |_{ [(2a+b)/3,(a+b)/2]}\cong \al _{A_n}$. 
\end{itemize}

\begin{lmm}\label{lmm:alpha}
Any basic loop $\be :[0,1]\to SC(Z)$ at $A_n$ is homotopic to the 
constant mapping $B$ in the subspace $p^{-1}([A, A_n]\times \mathbb{I})$.
\end{lmm}
\begin{proof}
We modify $\be$ to $\ga _0$ so that: 
\begin{enumerate}
\item[(1)] $\ga _0 (0)=\be (0)=A$, $\ga _0 |_{[1/3,1]} = \be |_{[1/3,1]}$; 
\item[(2)] $\ga _0 |_{[1/(4k+1), 1/(4k)]} \cong \ov{\al}_{A_{n+k}}$ and 
           $\ga _0 |_{[1/(4k+3), 1/(4k+2)]} \cong \al_{A_{n+k}}$ for
	   $k\ge 1$; and 
\item[(3)] $\ga _0 |_{ [1/4k,1/(4k-1)]}$ is a linear mapping and 
           $\ga _0 |_{ [1/(4k+2),1/(4k+1)]}$ is constant for $k\ge 1$. 
\end{enumerate}
It is easy to see that $\ga _0$ is homotopic to $\be$ in
 $p^{-1}([0,1/n]\times \{ 0\})$. This homotopy corresponds to the
 procedure I. Next we describe the homotopy corresponding to
 the procedure II according to the above classification
$(1)-(3)$. 
Let $E_{n,t}$ be the point $((t+n)/((n+1)n);t)$ on
 $L_{2n}$ and $F_{n,t}$ be the point $(1/n;t)$ on $L_{2n-1}$. 
We define $H:\mathbb{I}\times \mathbb{I}\to 
p^{-1}([0,1/n]\times \mathbb{I})$ so that $H(s,0) = \ga _0 (s)$ and the
 following hold:  
\smallskip

\item[(1)] $H(0,t)= (0;t)$, $H(s,t) = (2(1-s)/n; t)$ for
 $s\in[1/2,1]$ and $H(-,t)|_{[1/3,1/2]} \cong \al _{(1/n;t)}$; 
\item[(2)] $H(-,t)|_{[1/(4k+1), 1/4k]} \cong 
           \ov{\al} _{E_{n+k-1,t}}$ and $H(-,t)|_{[1/(4k+3), 1/(4k+2)]} 
           \cong \al_{F_{n+k,t}}$; and 
\item[(3)] $H(1/(4k+2),t) = F_{n+k,t}$, $H(1/(4k+1),t) = H(1/(4k),t) 
           = E_{n+k-1,t}$ and $H(-,t) |_{[1/4k,1/(4k-1)]}$ 
           and $H(-,t) |_{ [1/(4k+2),1/(4k+1)]}$ 
           are linear mappings. 

\smallskip
Then $H$ is continuous and it is a homotopy with $H(0,t) = H(1,t)$ for all
 $t$. Let $\ga _1 = H(-,1)$. Notice that 
$\ga _1|_{[1/(4k+1),1/(4k-2)]} \cong 
\ov{\al}_{B_{n+k-1}} c_{B_{n+k-1}}\al _{B_{n+k-1}}$ and 
$\ga _1|_{[1/(4k+2),1/(4k+1)]}$ is a linear mapping onto
 $[B_{n+k},B_{n+k-1}]$ and 
$\ga _1 |_{[1/2,1]}$ is a linear mapping onto $[B_n, B]$. 

It is then easy to see that $\ga _1$ is
 null-homotopic in $p^{-1}([0,1/n]\times \{ 1\})$, which corresponds to 
the procedure III. 
\end{proof}

%
%
\begin{center}
\beginpicture 
\setcoordinatesystem units <1.000mm,1.000mm>
\setplotarea x from -15 to 115, y from -15 to 115
\setsolid
\setplotsymbol ({\fiverm.})
\plot 11.604 98.220 8.930 50.008 /
\plot 11.315 0.012 11.604 98.220 /
\plot 13.840 2.365 13.390 100.000 /
\plot 20.267 100.207 13.840 2.365 /
\plot 23.403 96.840 16.061 0.000 /
\plot 22.890 0.044 23.403 96.840 /
\plot 27.088 4.176 26.580 100.000 /
\plot 41.612 100.802 27.705 4.131 /
\plot 46.306 94.439 31.279 -0.659 /
\plot 46.250 0.076 47.133 94.431 /
\plot 53.681 7.451 52.815 100.000 /
\plot 85.443 102.942 55.957 7.168 /
\plot 91.876 90.630 60.865 -2.372 /
\plot 86.670 0.000 94.164 90.016 /
\circulararc 171.313 degrees from 27.088 4.219 center at 27.399 4.199
\circulararc 180.000 degrees from 47.132 94.370 center at 46.719 94.405
\circulararc 177.246 degrees from 53.681 7.451 center at 54.822 7.337
\circulararc 159.523 degrees from 94.187 90.297 center at 93.002 90.255
\plot 9.830 100.000 8.040 50.032 /
\circulararc 179.831 degrees from 13.390 100.005 center at 11.610 100.000
\circulararc 175.791 degrees from 11.315 -0.000 center at 13.685 0.000
\circulararc 175.818 degrees from 26.580 100.023 center at 23.420 100.000
\circulararc 171.018 degrees from 22.890 -0.000 center at 27.110 0.000
\circulararc 171.277 degrees from 52.815 100.053 center at 47.185 100.000
\circulararc 162.140 degrees from 46.250 0.076 center at 53.750 0.000
\circulararc 158.145 degrees from 104.966 100.827 center at 95.000 100.000
\plot 113.330 0.000 104.966 100.827 /
\plot 100.000 0.000 100.000 100.000 /
\ellipticalarc axes ratio 4:1 360 degrees from 113.330 0.000 center at 100.000 0.000
\plot 0.000 0.000 100.000 0.000 /
\plot 100.000 100.000 0.000 100.000 /
\plot 0.000 100.000 0.000 0.000 /
\put {I} at 20 -10
\put {II} at -8 60
\put {III} at -8 110

\setdashes
\setplotsymbol ({\rm.})
\ellipticalarc axes ratio 4:1 360 degrees from 28.000 110.000 center at 13.000 110.000
\ellipticalarc axes ratio 1:4 360 degrees from 28.000 15.000 center at 28.000 0.000
\ellipticalarc axes ratio 4:1 360 degrees from 28.000 60.000 center at 13.000 60.000
\setsolid
\ellipticalarc axes ratio 4:1 360 degrees from 28.000 110.000 center at 13.000 110.000
\ellipticalarc axes ratio 1:4 270 degrees from 28.000 15.000 center at 28.000 0.000
\ellipticalarc axes ratio 4:1 -210 degrees from 28.000 60.000 center at 13.000 60.000

\endpicture
Figure 2.
\end{center}

{\it Proof of\/} Theorem~\ref{thm:main1}. 
By Lemma \ref{lmm:canonical}, we may start by a loop $f:[0,1]\to M_0$ with
 base point the $A$. Moreover, since the $A_m$'s are isolated points among
 $\{ A_m: m<\om \}$, are
 connected by intervals in $\mathbb{I}\times \{ 0\}$ and converge to
 $A$, every loop in $M_0$ is homotopic to a loop homotopic to an
 infinite concatenation of basic loops. 
We may assume that we have a disjoint family of intervals
 $(a_n,b_n)$ $(n<\nu )$, where $\nu \le \om$, such that each
 $f|_{[a_n,b_n]}$ is a basic loop at some $A_m$ and 
$\bigcup _{n<\nu} (a_n,b_n)$ is dense in $\mathbb{I}$. 
Then $f(s)=A$ for $s\notin \bigcup _{n<\nu} (a_n,b_n)$. 

We observe that the procedure II in the proof of Lemma~\ref{lmm:alpha} 
can be performed uniformly for all $f|_{[a_n,b_n]}$'s. Then we obtain a
 homotopy from $f$ to a loop in $M_1$ which consists of possibly
 infinitely many null-homotopic loops. Since the homotopies which
 correspond to the procedure III converge to $B$, we have a homotopy
 from $f$ to the constant mapping $B$. 
Extending the domain $[0,1]$ into both directions and adding a
 path from $A$ to $B$, we have a homotopy from $f$ to the constant
 mapping $A$ relative to their ends. 
\qed

\section{Proof of Theorem 1.2}

We shall show that $SC(Z)$ is noncontractible for every
noncontractible space $Z$. Let $SC_n(Z)$ be the
subspace of $SC(Z)$ defined as $p^{-1}(L_{2n-1}\cup
L_{2n}\cup L_{2n+1})$. 
\begin{dfn}
A mapping
$f:SC_n(Z)\to SC(Z)$ is said to be $\it{flat}$ if $p_2(f(z_1)) =
p_2(f(z_2))$, whenever $p_2(z_1) = p_2(z_2)$ for $z_1, z_2 \in
SC_n(Z)$. A homotopy $H: SC_n(Z)\times \mathbb{I} \to SC(Z)$ is
said to be flat if for every $t$, the mapping $H(-, t)$ is flat.
\end{dfn}
A similar idea of the flatness was used in \cite{KR}. 
\begin{lmm}\label{lmm:flat}
Let $n\in \mathbb{N}$ and let $H:SC_n(Z)\times\mathbb{I}\to SC(Z)$ be a
 mapping such that for every $y\in \mathbb{I}$ and $t\in \mathbb{I}$,
 the closure of the set
$p_2(H(M_y\cap SC_n(Z), t))$ does not contain both points
$0$ and $1$ and both mappings $H(-,0)$ and $H(-,1)$ are flat. 
Then there exists a flat homotopy from $H(-,0)$ to $H(-,1)$.
\end{lmm}

\begin{proof}
Fix the numbers $y$ and $t$. Let $A(y, t)$ and $B(y, t)$ be the
infimum and the supremum of the function $p_2\circ H(-, t):M_y\cap SC_n(Z)\to
\mathbb{I}$, respectively. Let
$$C(y,t) =\frac{A(y,t)}{1 + A(y,t) - B(y,t)}.$$ We note that $A(y,0) 
= C(y,0) = B(y,0)$ and $A(y,1) = C(y,1) = B(y,1)$, because $H(-,0)$ and
 $H(-,1)$ are flat. Consider the subset
$\mathbb{I}\times [A(y,t), B(y,t)].$ Let $\varphi$ be its
piecewise linear retraction to the interval $\mathbb{I}\times
\{C(y,t)\}$, which is defined by the mappings of vertices (see
Figure 1): $\varphi(B_{n,0}) = C_{2n-1}, \varphi(B_{n,1}) =
C_{2n}, \varphi(A_{n+1, 0}) = C_{2n}, \varphi(A_{n,1}) = C_{2n-1}$, 
where $A_{n,0}$ is the intersection point of $L_{2n-2}$ and 
$\{ (x; A(y,t)): x\in \mathbb{I}\}$, $A_{n,1}$ is the intersection point
 of $L_{2n-1}$ and $\{ (x; A(y,t)): x\in \mathbb{I}\}$, 
$B_{n,0}$ is the intersection point of $L_{2n-1}$ and $\{ (x; B(y,t))\}$ and 
$B_{n,1}$ is the intersection point of $L_{2n}$ and $\{ (x; B(y,t))\}$. 

If $A(y,t) = 0$ or $B(y, t) = 1,$ then $C(y, t) = 0$ or $C(y, t) =
1$ and $A_{n,0} = A_{n,1} = A_n$ or $B_{n,0} = B_{n,1} = B_n,$
respectively and the mapping $\varphi$ is well-defined.

Let $\psi_{y, t}$ be the natural retraction of
$p^{-1}(\mathbb{I}\times [A(y,t), B(y,t)])$ to
$p^{-1}(\mathbb{I}\times \{C(y,t)\})$ generated by
$\varphi$. Define now the homotopy $H':SC_n(Z)\times \mathbb{I}\to
SC(Z)$ by $H'(z,t) = \psi_{(p_2(z), t)}(H(z, t)).$ It is easy to
check that $p_2(H'(z, t)) = C(p_2(z), t)$ so that $H'$ is a flat
homotopy, and that $H'(-,0) = H(-,0)$ and $H'(-,1) = H(-,1)$. 
\end{proof}

To prove the Lemma~\ref{lmm:noncontr} below we introduce a notion which
will help us to
investigate flat homotopies. 

For $s\in (0,1)$ and $t\in \mathbb{I}$, we define a property
$P(s,t)$ of $H$ as follows: 
\begin{quote}
$H(M_s\cap SC_n(Z), t)\subseteq p^{-1}(\mathbb{I}\times (0,1))$ 
and the restriction of $H(-,t)$ to $M_s\cap SC_n(Z)$
is homotopic to the identity mapping on $M_s\cap SC_n(Z)$ in 
$p^{-1}(\mathbb{I}\times (0,1))$.
\end{quote}
We remark that by the flatness of $H$, if 
$H(M_s\cap SC_n(Z), t)\subseteq p^{-1}(\mathbb{I}\times (0,1))$,  
then there is a neighborhood $U$ of $(s;t)$ such that 
$H(M_{s'}\cap SC_n(Z), t')\subseteq p^{-1}(\mathbb{I}\times (0,1))$ for
any $(s';t')\in U$. 
\begin{lmm}\label{lmm:crucial}
Let $Z$ be a noncontractible space and $H:SC_n(Z)\times\mathbb{I}
\to SC(Z)$ a flat homotopy. If $p_2\circ H(M_0\cap SC_n(Z),t_0)
\subseteq p^{-1}(\mathbb{I}\times (0,1))$, then
 there exists a neighborhood $U$ of $(0,t_0)$ such that $H$ does not
 satisfy $P(s,t)$ for
 any $(s;t)\in U$ with $s>0$. An analogous statement holds for 
$p_2\circ H(M_1\cap SC_n(Z),t_0) \subseteq p^{-1}(\mathbb{I}\times (0,1))$. 
\end{lmm}
\begin{proof}
We have a neighborhood $U$ of $(0,t_0)$ such that 
$H(M_s\cap SC_n(Z), t)\subseteq p^{-1}(\mathbb{I}\times (0,1))$ for
any $(s;t)\in U$. We fix $(s;t)\in U$ with $s>0$ and assume to the contrary
 that $P(s,t)$ holds. 
Let $P_{2n+1},P_{2n},P_{2n-1}$ be the
 intersection of $\mathbb{I}\times \{ s\}$ and
 $L_{2n+1},L_{2n},L_{2n-1}$, respectively and 
$I_{2n} = [P_{2n+1},P_{2n}]$ and $I_{2n-1} = [P_{2n},P_{2n-1}]$. 
Then we have $M_s = \{ (0; s)\} \cup \bigcup _{n=1}^\infty 
I_n\cup Z_{P_n}$ and $M_s\cap SC_n(Z) = Z_{P_{2n+1}}\cup Z_{P_{2n}}
\cup Z_{P_{2n-1}}$. 
Since $H(-,t)$ maps $\bigcup _{u\in [P_{2n+1},A_{n+1}]\cup 
[A_{n+1},P_{2n}]}Z_u$ into $p^{-1}(\mathbb{I}\times (0,1))$, the restriction of $H(-,t)$ to 
$Z_{P_{2n+1}}\cup Z_{P_{2n}}$ is homotopic to a map
 $f:Z_{P_{2n+1}}\cup Z_{P_{2n}} \to \cup Z_{P_{2n}}$ 
in $p^{-1}(\mathbb{I}\times (0,1))$.

Since $M_s$ is a strong deformation retract of 
$p^{-1}(\mathbb{I}\times (0,1))$ similarly as in Lemma~\ref{lmm:retract}, 
$Z_{P_{2n+1}}\cup I_{2n} \cup Z_{P_{2n}}$ is a retract of 
$p^{-1}(\mathbb{I}\times (0,1))$. Since $Z$ is not contractible, the
 identity mapping on $Z_{P_{2n+1}} \cup Z_{P_{2n}}$ is not
 homotopic to any map  $f:Z_{P_{2n+1}} \cup Z_{P_{2n}} 
\to Z_{P_{2n}}$ in $p^{-1}(\mathbb{I}\times (0,1))$, which
 is a contradiction. 

To prove the statement for $H(M_1,t_0)$ we use $Z_{P_{2n}}\cup Z_{P_{2n-1}}$ 
and argue at a neighborhood of $B_n$ and obtain a similar conclusion.
\end{proof}
\begin{lmm}\label{lmm:noncontr}
Let $Z$ be a noncontractible space.
If $H:SC_n(Z)\times \mathbb{I}\to SC(Z)$ is a flat homotopy such
that $H(u, 0) = u$ for every $u\in SC_n(Z)$, then $H(-,1)$ is not a
 constant mapping.
\end{lmm}
\begin{proof}
To show this by contradiction, suppose that $H(-,1)$ is a constant
 mapping. 
Let $d: [0,1]\to S^1$ be a winding with the base point $s_0$, i.e. both
 $d |_{[0,1)}$ and $d |_{(0,1]}$ are bijective continuous mappings with
 $d(0)=d(1)=s_0$. 

We define a homotopy $H^{\ast}:S^1\times{\mathbb{I}}\to S^1$ as follows:
 \[
 H^{\ast}(u, t) = \left\{ 
\begin{array}{rl}
d(p_2\circ H(M_{d^{-1}(u)},t)) &\quad \mbox{ if }u\neq s_0 \mbox{ and }P(d^{-1}(u),t) \mbox{ holds;}\\
s_0, &\quad \mbox{ otherwise. }
\end{array}
\right.
\]
We have a contradiction with
 the fact that $H^{\ast}(s, 0) = s$ and that $S^1$ is not
 contractible, if $H^*$ is a homotopy (compare with \cite{Debski:Sine}). 
Hence it suffices to verify the continuity of $H^*$.

If $u\neq s_0$ and $P(d^{-1}(u),t)$ holds, the continuity at $(u,t)$ is
 clear. Otherwise, $u\neq s_0$ and $P(d^{-1}(u),t)$ does not hold, or 
$u = s_0$. We consider two cases: 
\noindent
Case 1.  Suppose that $u\neq s_0$ and $P(d^{-1}(u),t)$ does not hold: If 
$p_2\circ H(M_{d^{-1}(u)},t) = \{ 0\}$ or $\{ 1\}$, then the
 continuity at $(u,t)$ follows from that of $H$. Otherwise, since
 $H(-,t)$ maps $M_{d^{-1}(u)}\cap SC_n(Z)$ continuously with respect to
 $u$ and $t$, the restriction of $H(-,t)$ to $M_{d^{-1}(u)}\cap SC_n(Z)$
 is not homotopic to the identity on $M_{d^{-1}(u)}\cap SC_n(Z)$ in
 $p^{-1}(\mathbb{I}\times (0,1))$, i.e. $H^*$ takes the value $s_0$ in a
 neighborhood of $(u,t)$. 

\noindent
Case 2. Suppose that $u= s_0$: If each of $p_2\circ H(M_0\cap SC_n(Z),t)$ and
 $p_2\circ H(M_1\cap SC_n(Z),t)$ 
is equal to either $\{ 0\}$ or $\{ 1\}$, the continuity at $(u,t)$ follows
 from that of $H$. The remaining case is when 
$p_2\circ H(M_0\cap SC_n(Z),t)\subseteq (0,1)$ or 
$p_2\circ H(M_1\cap SC_n(Z),t)\subseteq (0,1)$. In this case the
 continuity follows by Lemma~\ref{lmm:crucial}. 
\end{proof}

{\it Proof of\/} Theorem~\ref{thm:main2}. 
To get a contradiction, suppose that $SC(Z)$ is a contractible space. Then
there exists a contraction $H:SC(Z)\times \mathbb{I}\to SC(Z)$. By the
compactness of the time interval $\mathbb{I}$, for every $a = (0;y) \in \{ 0\}
\times \mathbb{I}$, there exists $\vep _0>0$ such that the diameters of
$p_2\circ H(O_{\vep _0}(a),t)$ are less than $1$ for all $t\in
\mathbb{I}$. Hence, by compactness of $L$, 
there exists $\vep _1> 0$ such that the diameters 
$p_2\circ H(O_{\vep _1}(a),t)$ are less than $1$ for all $a = (0;y)\in \{ 0\}
\times \mathbb{I}$ and all $t\in \mathbb{I}$. 
Let $n$ be a number such that $1/n < \vep _1$. By Lemma \ref{lmm:flat}
we may then assume that $H|_{SC_n(Z)\times \mathbb{I}}$ is a flat
contraction. However, this contradicts Lemma \ref{lmm:noncontr}.
\qed

\begin{rmk}
The space $SC(S^1)$ is simply connected and it follows
by the Mayer-Vietoris exact sequence for the singular homology 
 that $H_n(SC(S^1)) = 0$ for $n\ge 3$. 
The question whether $H_2(SC(S^1)) = 0$ was open in the first
 draft of this paper. This has subsequently been answered negatively by the
 authors and also independently by J. Dydak and A. Mitra. The proof will
 appear in our forth-coming paper ``An example of a nonaspherical cell-like
2-dimensional continuum and some related constructions''. 
\end{rmk}

\section{Acknowledgements}

The authors thank the referee for careful reading, detecting errors,
and suggestions to make arguments clearer. 
They also thank Haruto Ohta for discussions about dimension theory. 
The first and the third author were supported by the
Japanese-Slovenian research grant BI-JP/03-04/2 and 
the Grant-in-Aid for Scientific research (C) of Japan 
No. 16540125. 

\vfill\eject

\end{document}